\documentclass[]{article}

\usepackage{amssymb,amsmath,hyperref}
\usepackage{color,graphics,graphicx}
\usepackage{color}
\usepackage{amsthm}	
\usepackage[margin=1in]{geometry}

\title{Robust $H_{\infty}$ Loop-Shaping Differential Thrust Control Methodology for Lateral/Directional Stability of an Aircraft with a Damaged Vertical Stabilizer}

\author{
  Long Lu%
    \thanks{Graduate Student, Department of Aerospace Engineering, AIAA Student Member, long.lu@sjsu.edu.}
\  and \ Kamran Turkoglu\thanks{Assistant Professor, Department of Aerospace Engineering, AIAA Member, kamran.turkoglu@sjsu.edu}\\
  {\normalsize\itshape San Jos\'{e} State University, San Jose, CA 95192, USA} \\ \\
 }

\begin{document}

\maketitle

\begin{abstract}

The vertical stabilizer is the key aerodynamic surface that provides an aircraft with its directional stability characteristic while ailerons and rudder are the primary control surfaces that give pilots control authority of the yawing and banking maneuvers. Losing the vertical stabilizer will, therefore, result in the consequential loss of lateral/directional stability and control, which is likely to cause a fatal crash. In this paper, we construct a scenario of a damaged aircraft model which has no physical rudder control surface, and then a strategy based on differential thrust is proposed to be utilized as a control input to act as a "virtual" rudder to help maintain stability and control of the damaged aircraft. The $H_{\infty}$ loop-shaping approach based robust control system design is implemented to achieve a stable and robust flight envelope, which is aimed to provide a safe landing. Investigation results demonstrate successful application of such robust differential thrust methodology as the damaged aircraft can achieve stability within feasible control limits. Finally, the robustness analysis results conclude that the stability and performance of the damaged aircraft in the presence of uncertainty remain within desirable limits, and demonstrate not only a robust, but a safe flight mission through the proposed $H_{\infty}$ loop-shaping robust differential thrust control methodology.

\end{abstract}

\section*{Nomenclature}

\begin{tabbing}
  XXX \= \kill
$A/C$ 		\> aircraft \\
$A_d$ 		\> state matrix of the damaged aircraft \\
$B_d$ 		\> input matrix of the damaged aircraft \\
$b$ 		\> aircraft wing span, $ft$ \\ 
$C_d$ 		\> output matrix of the damaged aircraft\\
${C_L}_{i}$		\> dimensionless derivative of rolling moment, $i=p, r, \beta, \delta_a, \delta_r$\\
${C_N}_{i}$		\> dimensionless derivative of yawing moment, $i=p, r, \beta, \delta_a, \delta_r$\\
${C_Y}_{i}$		\> dimensionless derivative of side force, $i=p, r, \beta, \delta_a, \delta_r$\\
$D_d$ 		\> state transition matrix of the damaged aircraft \\
$\frac{d\sigma}{d\beta}$ 		\> change in side wash angle with respect to change in side slip angle\\
$g$ 		\> gravitational acceleration, $ft/s^2$ \\
$\overline{I_{xx}}$ 		\> normalized mass moment of inertia about the x axis, $slug*ft^2$ \\
$\overline{I_{xz}}$ 		\> normalized product of inertia about the xz axis, $slug*ft^2$\\
$\overline{I_{zz}}$ 		\> normalized mass moment of inertia about the z axis, $slug*ft^2$ \\
$L_{i}$		\> dimensional derivative of rolling moment, $i=p, r, \beta, \delta_a, \delta_r$\\
$L_v$ 		\>  vertical stabilizer lift force, $lbf$\\
$l_v$ 		\> distance from the vertical stabilizer aerodynamic center to the aircraft center of gravity, $ft$\\
$m$ 		\> aircraft mass, $slugs$ \\
$N_{i}$		\> dimensional derivative of yawing moment, $i=p, r, \beta, \delta_a, \delta_r, \delta T$\\
$p$ 		\> roll rate, $deg/s$ \\
$r$ 		\> yaw rate, $deg/s$ \\
$S$ 		\> aircraft wing area, $ft^2$ \\
$S_v$ 		\> vertical stabilizer area, $ft^2$ \\
$T$ 		\> engine thrust, $lbf$ \\
$T_c$ 		\> engine thrust command, $lbf$ \\
$t$		\> time, $s$ \\
$t_d$		\> time delay, $s$\\
$V_v$ 		\> vertical stabilizer volume ratio, $ft^3$\\
$v$ 		\> airspeed, $ft/s$ \\
$W$ 		\> aircraft weight, $lbs$ \\
$Y_{i}$		\> dimensional derivative of side force, $i=p, r, \beta, \delta_a, \delta_r$\\
$y_e$ 		\> distance from the outermost engine to the aircraft center of gravity, $ft$\\
$z_v$ 		\> distance from the vertical stabilizer center of pressure to the fuselage center line, $ft$\\
$\alpha$ 		\> angle of attack, $deg$ \\
$\beta$ 		\> side slip angle, $deg$ \\
$\gamma$  	\> flight path angle, $deg$ \\
$\delta_a$		\> aileron deflection, $deg$ \\
$\delta_r$		\> rudder deflection, $deg$ \\
$\Delta T$		\> collective thrust, $lbf$ \\
$\delta T$		\> differential thrust, $lbf$ \\
$\zeta$		\> damping ratio \\
$\eta$  	\> efficiency factor \\
$\theta$  	\> pitch angle, $deg$ \\
$\rho$		\> air density, $slug/ft^3$ \\
$\tau$   	\> time constant, $s$ \\
$\phi$    	\> roll angle, $deg$ \\
$\omega$		\> bandwidth frequency , $1/s$ \\
$\dot{( \ )}$ 	\>  first order time derivative \\
$\ddot{( \ )}$ 	\>  second order time derivative \\
$\bar{( \ )}$ 	\>  trimmed value \\
${( \ )}_d$ 	\>  damaged aircraft dynamics \\
 \end{tabbing}

\section{Introduction} \label{Introduction} 

The vertical stabilizer of an aircraft is an essential element in providing the aircraft with its directional stability characteristic while ailerons and rudder serve as the primary control surfaces of the yawing and banking maneuvers. In the event of an aircraft losing its vertical stabilizer, the sustained damage will cause lateral/directional stability to be compromised and lack of control is likely to result in a fatal crash. Notable examples of such a scenario are the crash of the American Airline 587 in 2001 when an Airbus A300-600 lost its vertical stabilizer in wake turbulence, killing all passengers and crew members \cite{NTSB_flight587} and the crash of Japan Airlines Flight 123 in 1985 when a Boeing 747-SR100 lost its vertical stabilizer leading to an uncontrollable aircraft, resulting in 520 casualties \cite{faa_flight123}.

%
 
However, not all situations of losing the vertical stabilizer resulted in a total disaster. In one of those cases, the United Airlines Flight 232 in 1989 \cite{NTSB_flight_232}, differential thrust was proved to be able to make the aircraft controllable. Another remarkable endeavor is the landing of the Boeing 52-H even though the aircraft lost most of its vertical stabilizer in 1964 \cite{Hartnett06}.


In literature, research on this topic has been conducted with two main goals: to understand the response characteristics of the damaged aircraft such as the work of Bacon and Gregory \cite{Bacon_Gregory_07}, Nguyen and Stepanyan \cite{NguyenStepanyan10}, and Shah \cite{Shah08}, as well as to come up with an automatic control algorithm to save the aircraft from disasters where the work of Burcham et al. \cite{Burcham_et_al_98}, Guo et al. \cite{Guo_et_al_11}, Liu et al. \cite{Liu_et_al_09}, Tao and Ioanou \cite{Tao_Ioannou_88}, and Urnes and Nielsen \cite{UrnesNielsen10} hold the detailed analysis.

There are many valuable efforts in literature, such as the work of Shah \cite{Shah08}, where a wind tunnel study was performed to evaluate the aerodynamic effects of damages to lifting and stability/control surfaces of a commercial transport aircraft. In his work, Shah \cite{Shah08} studied this phenomenon in the form of partial or total loss of the wing, horizontal, or vertical stabilizers for the development of flight control systems to recover the damaged aircraft from adverse events. The work of Nguyen and Stepanyan \cite{NguyenStepanyan10} investigates the effect of the engine response time requirements of a generic transport aircraft in severe damage situations associated with the vertical stabilizer. They carried out a study which concentrates on assessing the requirements for engine design for fast responses in an emergency situation. In addition, the use of differential thrust as a propulsion command for the control of directional stability of a damaged transport aircraft was studied by Urnes and Nielsen \cite{UrnesNielsen10} to identify the effect of the change in aircraft performance due to the loss of the vertical stabilizer and to make an improvement in stability utilizing engine thrust as an emergency yaw control mode with feedback from the aircraft motion sensors.  

Existing valuable research in literature provides very unique insight into the dynamics of such an extreme scenario whereas in this paper a unique extension to the existing works is provided. Here, we provide a $H_{\infty}$ loop-shaping based robust differential thrust control methodology to control such a damaged aircraft and land safely. This research is motivated to improve air travel safety in commercial aviation by incorporating the utilization of differential thrust to regain lateral/directional stability for a commercial aircraft (in this case, a Boeing 747-100) with a damaged vertical stabilizer. For this purpose, the damaged aircraft model is constructed in Section \ref{AC Model}, where its plant dynamics is investigated. In Section \ref{Differential Thrust}, the engine dynamics of the aircraft is modeled as a system of differential equations with corresponding time constant and time delay terms to study the engine response characteristic with respect to a differential thrust input, and the novel differential thrust control module is developed to map the rudder input to differential thrust input. In Section \ref{Open Loop}, the aircraft's open loop system response is investigated. In Section \ref{Robust}, the robust control system design based on $H_{\infty}$ loop-shaping approach is implemented as a means to stabilize the damaged aircraft. Investigation provides remarkable results of such robust differential thrust methodology as the damaged aircraft can reach steady state stability relatively fast, under feasible control efforts. In Section \ref{Robustness}, robustness and uncertainty analysis is conducted to test the guaranteed stability bounds and to validate the overall performance of the system in the presence of given ball of uncertainty. With Section \ref{Conclusion}, the paper is finalized, and future work is discussed. 

\section{The Damaged Aircraft Model} \label{AC Model}
In this paper, the flight scenario is chosen to be a steady, level cruise flight for the Boeing 747-100 at Mach 0.65 and 20,000 feet. It is assumed that at one point during the flight, the vertical stabilizer is completely damaged, and in the followings the means to control such aircraft is investigated in such an extreme case scenario. For this purpose, first, the damaged aircraft model is developed for analysis, where the flight conditions for the damaged model are summarized in Table \ref{table:Flight conditions}.
\begin{table}[htbp!]
\centering
\caption{Flight conditions}
\label{table:Flight conditions}
 \begin{tabular}{|c|c|}
 \hline
 \textbf{Parameter} & \textbf{Value}\\ \hline
 Altitude $(ft)$ & 20,000 \\ \hline
Air Density $(slug/ft^3)$ & 0.001268 \\ \hline
Airspeed $(ft/s)$ & 673 \\ \hline
\end{tabular}
\end{table}

\subsection{The Damaged Aircraft Model} \label{AC model}

The Boeing 747-100 was chosen for this research due to its widely available technical specification, aerodynamics, and stability derivative data. The data for the nominal (undamaged) Boeing 747-100 are summarized in Table \ref{table:The undamaged aircraft data}.

\begin{table}[htbp!]
\centering
\caption{The nominal (undamaged) aircraft data \cite{747_ac_charac, roskam87}}
\label{table:The undamaged aircraft data}
 \begin{tabular}{|c|c|}
 \hline
 \textbf{Parameter} & \textbf{Value}\\ \hline
$S$ $(ft^2)$ & 5500  \\ \hline
$b$ $(ft)$ & 196  \\ \hline
$\bar{c}$ $(ft)$ & 27.3  \\ \hline
$y_e$ $(ft)$ & 69.83  \\ \hline
$W$ $(lbs)$ &  $6.3663*10^5$\\ \hline
$m$ $(slugs)$ &  $19786.46$\\ \hline 
$I_{xx}$ $(slug*ft^2)$ & $18.2*10^6$  \\ \hline 
$I_{yy}$ $(slug*ft^2)$ & $33.1*10^6$  \\ \hline
$I_{zz}$ $(slug*ft^2)$ & $49.7*10^6$  \\ \hline
$I_{xz}$ $(slug*ft^2)$ & $0.97*10^6$  \\ \hline
${C_L}_{\beta}$	& -0.160 \\ \hline 
${C_L}_{p}$	& -0.340 \\ \hline 
${C_L}_{r}$	& 0.130 \\ \hline 
${C_L}_{\delta_a}$	& 0.013 \\ \hline 
${C_L}_{\delta_r}$	& 0.008 \\ \hline 
${C_N}_{\beta}$	& 0.160 \\ \hline 
${C_N}_{p}$	& -0.026 \\ \hline 
${C_N}_{r}$	& -0.28 \\ \hline 
${C_N}_{\delta_a}$	& 0.0018 \\ \hline 
${C_N}_{\delta_r}$	& -0.100 \\ \hline 
${C_Y}_{\beta}$	& -0.90 \\ \hline 
${C_Y}_{p}$	& 0 \\ \hline 
${C_Y}_{r}$	& 0 \\ \hline 
${C_Y}_{\delta_a}$	& 0 \\ \hline 
${C_Y}_{\delta_r}$	& 0.120 \\ \hline 
\end{tabular}
\end{table} 

For the modeling of the damaged aircraft, in case of the loss of the vertical stabilizer, lateral/directional stability derivatives need to be re-examined and recalculated due to the resulting significant changes in the geometry and the aerodynamic behavior. Since the whole aerodynamic structure is affected, the new corresponding stability derivatives have to be calculated and studied. The lateral/directional dimensionless derivatives that depend on the vertical stabilizer include \cite{nelson98}:

\begin{equation}
 {C_Y}_{\beta}=-\eta \frac{S_v}{S}{{C_L}_\alpha}_v\left(1+\frac{d\sigma}{d\beta}\right)
\end{equation}
 
\begin{equation}
  {C_Y}_r=-2\left(\frac{l_v}{b}\right){{C_Y}_\beta}_{tail}
\end{equation}
   
\begin{equation}   
 {C_N}_{\beta}={{C_N}_{\beta}}_{wf}+\eta_v V_v{{C_L}_\alpha}_v\left(1+\frac{d\sigma}{d\beta}\right)
\end{equation}

\begin{equation}
  {C_N}_r=-2\eta_v V_v\left(\frac{l_v}{b}\right){{C_L}_{\alpha}}_v
\end{equation}

\begin{equation}
  {C_L}_r=\frac{C_L}{4}-2\left(\frac{l_v}{b}\right)\left(\frac{z_v}{b}\right){{C_Y}_\beta}_{tail}
\end{equation}

When the vertical stabilizer is completely lost in an aircraft, the vertical tail area, vertical tail volume, and associated efficiency factor will all be zero, resulting in ${C_Y}_{\beta}={C_Y}_r={C_N}_r=0$. If the vertical stabilizer is assumed to be the primary aerodynamic surface responsible for the weathercock stability, then ${C_N}_{\beta}=0$ and finally, we have ${C_L}_r=\frac{C_L}{4}$.

In addition, without the vertical stabilizer, the mass and inertia data of the damaged aircraft are going to change as well, where the values that reflect such a scenario (for the damaged aircraft) are calculated and listed in Table \ref{table:The damaged aircraft mass and inertia data}.

\begin{table}[htbp!]
\centering
\caption{The damaged aircraft mass and inertia data}
\label{table:The damaged aircraft mass and inertia data}
 \begin{tabular}{|c|c|}
 \hline
 \textbf{Parameter} & \textbf{Value}\\ \hline
$W$ $(lbs)$ &  $6.2954*10^5$\\ \hline 
$m$ $(slugs)$ &  $19566.10$\\ \hline 
$I_{xx}$ $(slug*ft^2)$ & $17.893*10^6$  \\ \hline 
$I_{yy}$ $(slug*ft^2)$ & $30.925*10^6$  \\ \hline
$I_{zz}$ $(slug*ft^2)$ & $47.352*10^6$  \\ \hline
$I_{xz}$ $(slug*ft^2)$ & $0.3736*10^6$  \\ \hline
\end{tabular}
\end{table} 

In this study, during an event of the loss of the vertical stabilizer, it is proposed that the differential thrust component of aircraft dynamics be utilized as an alternate control input replacing the rudder control to regain stability and control of lateral/directional flight dynamics. Next, the lateral-directional linear equations of motion of the damaged aircraft are presented, with the ailerons $(\delta_a)$, differential thrust $(\delta T)$, and collective thrust $(\Delta T)$ as control inputs \cite{NguyenStepanyan10}, as 

\begin{equation}
\left[
\begin{array}{c}
\dot{\phi}\\
\dot{p}\\
\dot{\beta}\\
\dot{r}
\end{array}\right]
= \left[
\begin{array}{cccc}
0 & 1 & 0 & \bar{\theta} \\
0 & L_p & L_{\beta} & L_r \\
\frac{g}{\bar{V}} & \frac{Y_p}{\bar{V}} & \frac{Y_{\beta}+g\bar{\gamma}}{\bar{V}} & \frac{Y_p}{\bar{V}}-1 \\
0 & N_p & N_{\beta} & N_r
\end{array}\right]
\left[
\begin{array}{c}
{\phi}\\
{p}\\
{\beta}\\
{r}
\end{array}\right]
+
\left[
\begin{array}{ccc}
0 & 0 & 0\\
L_{\delta_a} & \frac{\overline{I_{xz}}y_e}{\overline{I_{xx}}\overline{I_{zz}}-\overline{I_{xz}}^2} & 0\\
\frac{Y_{\delta_a}}{\bar{V}} & 0 & \frac{-\bar{\beta}}{m\bar{V}}\\
N_{\delta_a} & \frac{\overline{I_{xx}}y_e}{\overline{I_{xx}}\overline{I_{zz}}-\overline{I_{xz}}^2} & 0
\end{array}\right]
\left[
\begin{array}{c}
{\delta_a}\\
{\delta T} \\
{\Delta T} \\
\end{array}\right]
\end{equation}

In this case, if the initial trim side-slip angle is zero, then $\Delta T$ does not have any significance in the control effectiveness for a small perturbation around the trim condition \cite{NguyenStepanyan10}, which means that the above equations can be reduced to the final form as:

\begin{equation}
\left[
\begin{array}{c}
\dot{\phi}\\
\dot{p}\\
\dot{\beta}\\
\dot{r}
\end{array}\right]
= \left[
\begin{array}{cccc}
0 & 1 & 0 & \bar{\theta} \\
0 & L_p & L_{\beta} & L_r \\
\frac{g}{\bar{V}} & \frac{Y_p}{\bar{V}} & \frac{Y_{\beta}+g\bar{\gamma}}{\bar{V}} & \frac{Y_p}{\bar{V}}-1 \\
0 & N_p & N_{\beta} & N_r
\end{array}\right]
\left[
\begin{array}{c}
{\phi}\\
{p}\\
{\beta}\\
{r}
\end{array}\right]
+
\left[
\begin{array}{cc}
0 & 0\\
L_{\delta_a} & \frac{\overline{I_{xz}}y_e}{\overline{I_{xx}}\overline{I_{zz}}-\overline{I_{xz}}^2}\\
\frac{Y_{\delta_a}}{\bar{V}} & 0 \\
N_{\delta_a} & \frac{\overline{I_{xx}}y_e}{\overline{I_{xx}}\overline{I_{zz}}-\overline{I_{xz}}^2}
\end{array}\right]
\left[
\begin{array}{c}
{\delta_a}\\
{\delta T} \\
\end{array}\right]
\end{equation}

\subsection{Plant Dynamics} \label{Plant Dynamics}

%
%
%
%

Based on the calculated data for the lateral/directional stability derivatives, mass and inertia of the aircraft without its vertical stabilizer, the lateral/directional state-space representation of the damaged aircraft is obtained as the following

\begin{equation} \label{eq:Ad}
{A_d}= \left[
\begin{array}{cccc}
0 & 1 & 0 & 0\\
0 & -0.8566 & -2.7681 &  0.1008\\
0.0478 & 0 & 0 & -1 \\
0 & -0.0248 & 0 & 0
\end{array}\right]
\end{equation}

\begin{equation} \label{eq:Bd}
{B_d}= \left[
\begin{array}{cc}
0 & 0\\
0.2249 & 0.0142\\
0 & 0\\
0.0118 & 0.6784
\end{array}\right]
\end{equation}

\begin{equation} \label{eq:Cd}
C_d= \left[
\begin{array}{cccc}
1 & 0 & 0 & 0\\
0 & 1 & 0 & 0\\
0 & 0 & 1 & 0\\
0 & 0 & 0& 1
\end{array}\right]
\end{equation}

\begin{equation} \label{eq:Dd}
D_d= \left[
\begin{array}{cc}
0 & 0\\
0 & 0\\
0 & 0\\
0 & 0
\end{array}\right]
\end{equation}
where the subscript $(~~)_{d}$ stands for the $damaged$ aircraft dynamics.

Here, ${A_d}$ represents the state matrix of the damaged aircraft, and ${B_d}$ stands for the input matrix of the scenario where the ailerons $(\delta_a)$ and differential thrust $(\delta T)$ are control inputs of the damaged aircraft. ${C_d}$ is the output matrix, and ${D_d}$ is the direct (state transition) matrix of the damaged aircraft.

In addition, for more in-depth intuition, the damping characteristics of the damaged aircraft are also provided in Table \ref{table:damaged_damping}.


\begin{table}[htbp!]
\centering
\caption{Damping characteristics of the damaged aircraft}
\label{table:damaged_damping}
 \begin{tabular}{|c|c|c|c|c|}
 \hline
\textbf{Mode} & \textbf{Pole Location} & \textbf{Damping} & \textbf{Frequency $(1/s)$} & \textbf{Period $(s)$}\\ \hline 
Dutch Roll & $0.0917\pm i0.43$ & $-0.209$ & $0.439$ & $14.2969$\\ \hline
Spiral & $6.32*10^{-18}$ & $-1$ & $6.32*10^{-18}$ & $9.9486*10^{17}$\\ \hline
Roll & $-1.04$ & $1$ & $1.04$ & $6.0422$\\
\hline
\end{tabular}
\end{table}

Table \ref{table:damaged_damping} clearly indicates the unstable nature of the damaged aircraft in the Dutch roll mode by the obvious Right Half Plane (RHP) pole locations. Furthermore, the pole of the spiral mode lies at the origin. The only stable mode of the damaged aircraft remains to be the roll mode. 

Following the aerodynamic and stability analysis of the damaged aircraft, we next concentrate on the differential thrust methodology.



\section{Differential Thrust as a Control Mechanism} \label{Differential Thrust}

\subsection{Propulsion Dynamics} \label{Propulsion}

With emerging advancements in manufacturing processes, structures, and materials, it is a well known fact that aircraft engines have become highly complex systems and include numerous nonlinear processes, which affect the overall performance (and stability) of the aircraft. From the force-balance point of view, this is usually due to the existing coupled and complex dynamics between engine components and their relationships in generating thrust. However, in order to utilize the differential thrust generated by the jet engines as a control input for lateral/directional stability, the dynamics of the engine need to be modeled in order to gain an insight into the response characteristics of the engines.  

Engine response, generally speaking, depends on its time constant and time delay characteristics. Time constant dictates how fast the thrust is generated by the engine while time delay (which is inversely proportional to the initial thrust level) is due to the lag in engine fluid transport and the inertia of the mechanical systems such as rotors and turbo-machinery blades \cite{NguyenStepanyan10}. 

It is also suggested \cite{NguyenStepanyan10} that the non-linear engine dynamics model can be simplified as a second-order time-delayed linear model as
\begin{equation}
\ddot{T}+2\zeta\omega\dot{T}+\omega^2T= \omega^2T_c(t-t_d)
\end{equation}
where $\zeta$ and $\omega$ are the damping ratio and bandwidth frequency of the closed-loop engine dynamics, respectively; $t_d$ is the time delay factor, and $T_c$ is the thrust command prescribed by the engine throttle resolver angle.

With the time constant defined as the inverse of the bandwidth frequency $(\tau=\frac{1}{\omega})$, and $\zeta$ chosen to be 1 representing a critically damped engine response (to be comparable to existing studies), the engine dynamics can be represented as

\begin{equation}
\left[
\begin{array}{c}
\dot{T}\\
\ddot{T}
\end{array}\right]
= \left[
\begin{array}{cc}
0 & 1\\
\frac{-1}{\tau^2} & \frac{-2}{\tau}
\end{array}\right]
\left[
\begin{array}{c}
T\\
\dot{T}
\end{array}\right]
+
\left[
\begin{array}{c}
0 \\
\frac{1}{\tau^2}
\end{array}\right]
T_c(t-t_d)
\end{equation}

For this study, the Pratt and Whitney JT9D-7A engine is chosen for the application in the Boeing 747-100, where the engine itself produces a maximum thrust of 46,500 lbf \cite{boeing747_tech_spec}. At Mach 0.65 and 20,000 feet flight conditions, the engine time constant is 1.25 seconds, and the time delay is 0.4 second \cite{NguyenStepanyan10}.

\begin{figure}[h!]
 \centering
 \includegraphics[width=4truein]{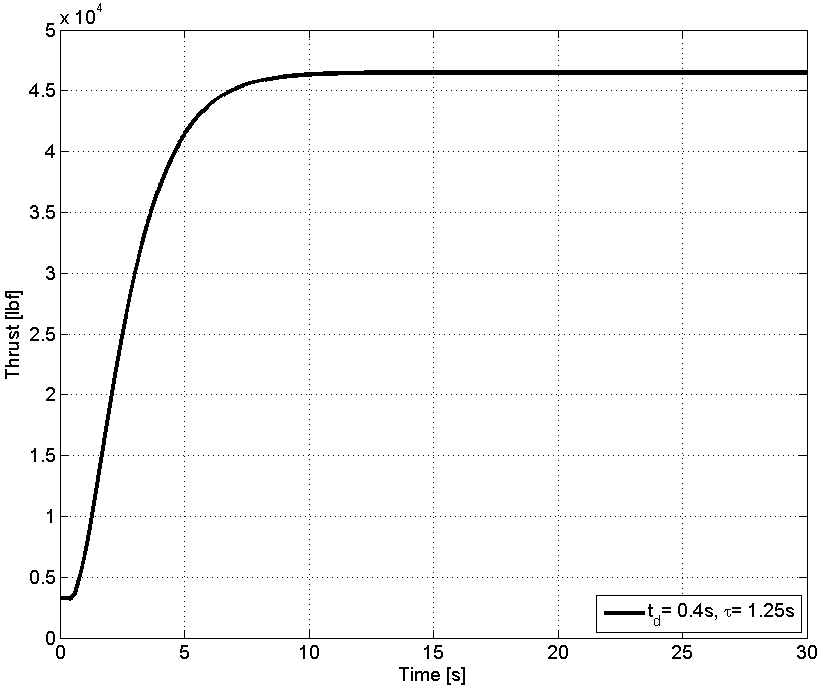}
 \caption{Engine thrust response at Mach 0.65 and 20,000 feet}
 \label{fig:engine_response}
 \end{figure}

The engine thrust response curve at Mach 0.65 and 20,000 feet is, therefore, obtained as shown in Fig. \ref{fig:engine_response}, which provides a useful insight into how the time constant and time delay factors affect the generation of thrust for the JT9D-7A jet engine. At Mach 0.65 and 20,000 feet, with the engine time constant of 1.25 seconds, and the time delay of 0.4 second, it takes approximately ten seconds for the engine to reach steady state and generate its maximum thrust capacity at 46,500 lbf from the trim thrust of 3221 lbf. The increase in thrust generation follows a relatively linear fashion with the engine response characteristic of approximately 12,726 lbf/s during the first two seconds, and then the thrust curve becomes nonlinear until it reaches its steady state at maximum thrust capacity after about ten seconds. This represents one major difference between the rudder and differential thrust as a control input. Due to the lag in engine fluid transport and turbo-machinery inertia, differential thrust (as a control input) cannot respond as instantaneously as the rudder, which has to be taken into account very seriously in control system design.

\subsection{Thrust Dynamics and Configuration}
\begin{figure}[htbp!]
 \centering
 \includegraphics[width=6truein]{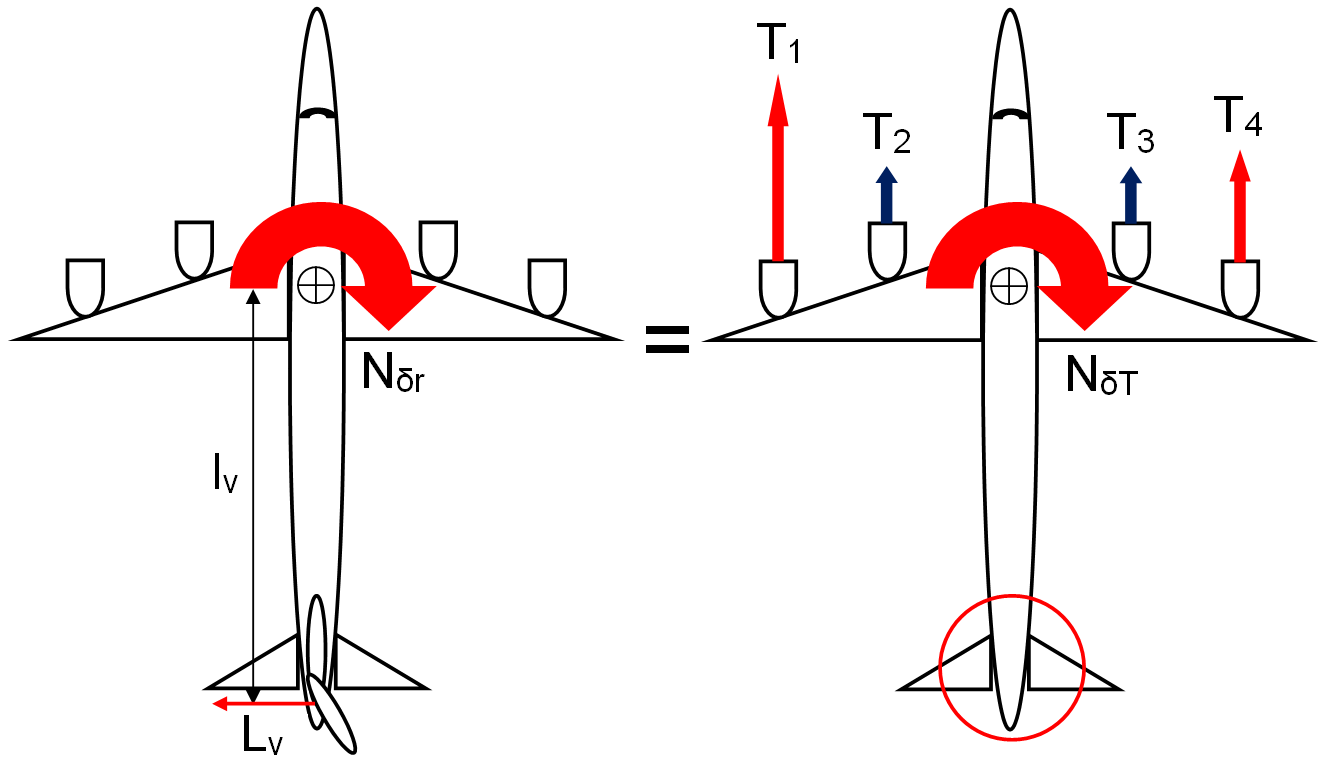}
 \caption{The free body diagram}
 \label{fig:model}
 \end{figure}

In order to utilize differential thrust as a control input for the four-engined Boeing 747-100 aircraft, a differential thrust control module must be developed. Here, the differential thrust input is defined as the difference between the thrust generated by engine number 1 and engine number 4 while the amounts of thrust generated by engine number 2 and 3 are kept equal to each other as shown in Eqs. (\ref{eq:d_th}-\ref{eq:th23}). This concept is illustrated in further details in Fig. \ref{fig:model}.

\begin{equation}\label{eq:d_th}
\delta T=T_1-T_4
\end{equation}
\begin{equation}\label{eq:th23}
T_2=T_3
\end{equation}

Engine number 1 and 4 are employed to generate the differential thrust due to the longer moment arm $(y_e)$, which makes the differential thrust more effective as a control for yawing moment. This brings into the picture the need of developing a logic that maps the rudder input to differential thrust input, which is further explained in the following section.

\subsection{Rudder Input to Differential Thrust Input Mapping Logic}

\begin{figure}[htbp!]
 \centering
 \includegraphics[width=5truein]{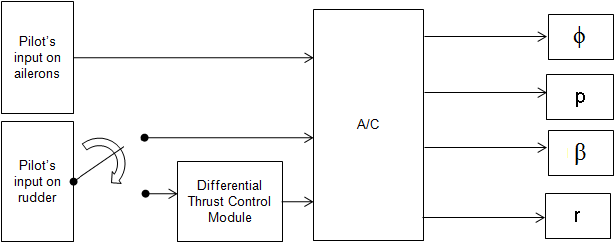}
 \caption{Aircraft control logic diagram}
 \label{fig:control_logics}
 \end{figure}

\begin{figure}[htbp]
 \centering
 \includegraphics[width=5truein]{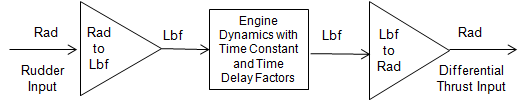}
 \caption{Differential thrust control module}
 \label{fig:dT_control_module}
 \end{figure}

When the vertical stabilizer of the aircraft is intact (i.e. with nominal plant dynamics), the pilot has the ailerons and rudder as major control surfaces. However, when the vertical stabilizer is damaged, most probably, the pilot will keep on demanding the control effort from the rudder until it is clear that there is no response from the rudder. To eliminate this mishap, but to still be able to use the rudder demand, a differential thrust control module is introduced in the control logic, as shown in Fig. \ref{fig:control_logics} and Fig. \ref{fig:dT_control_module}, respectively. This differential thrust control module maps  corresponding input/output dynamics from the rudder pedals to the aircraft response, so that when the rudder is lost, the rudder input (from the pilot) will still be utilized but switched to the differential thrust input, which will act as a \emph{rudder-like input} for lateral/directional controls. This logic constitutes one of the novel approaches introduced in this paper.

As it can be also seen from Fig. \ref{fig:control_logics} and Fig. \ref{fig:dT_control_module}, the differential thrust control module's function is to convert the rudder pedal input from the pilot to the differential thrust input. In order to achieve that, the rudder pedal input (in radians) is converted to the differential thrust input (in pounds-force) which is then provided into the engine dynamics, as discussed previously in Section \ref{Propulsion}. With this modification, the engine dynamics will dictate how differential thrust is generated, which is then provided as a "virtual rudder" input into the aircraft dynamics. The radian to pound-force mapping is provided in the next section, for the completeness of analysis.

\subsection{Radian to Pound-Force Conversion Factor}
Using Fig. \ref{fig:model} and with the steady, level flight assumption, the following relationship can be obtained:
\begin{equation}
N_{\delta_r}= N_{\delta T}
\end{equation}

\begin{equation}
qSb{C_N}_{\delta_r}\delta_r=(\delta T) y_e
\end{equation}
which means the yawing moment by deflecting the rudder and by using differential thrust have to be the same. Therefore, the relationship between the differential thrust control input $(\delta T)$ and the rudder control input $(\delta_r)$ can be obtained as

\begin{equation}
\delta T= \left(\frac{qSb{C_N}_{\delta_r}}{y_e}\right)\delta_r
\end{equation}

Based on the flight conditions at Mach 0.65 and 20,000 feet, and the data for the Boeing 747-100 summarized in Table \ref{table:Flight conditions} and Table \ref{table:The undamaged aircraft data}, the conversion factor for the rudder control input to the differential thrust input is calculated to be 

\begin{equation}
\frac{\delta T}{\delta_r}= -4.43*10^5 \> \frac{lbf}{rad}
\end{equation}

Due to the sign convention of rudder deflection and the free body diagram in Fig. \ref{fig:model}, $\delta_r$ here is negative. Therefore, for the Boeing 747-100, in this study, the conversion factor for the mapping of a rudder input to a differential thrust input is found to be

\begin{equation}
\frac{\delta T}{\delta_r}= 4.43*10^5 \> \frac{lbf}{rad}
\end{equation}

\subsection{Commanded vs. Available Differential Thrust}
At this point, the worst case scenario is considered, and it is assumed that the aircraft has lost its vertical stabilizer so that the rudder input is converted to the differential thrust input according to the logic discussed previously in this section. 

\begin{figure}[htbp!]
 \centering
 \includegraphics[width=4truein]{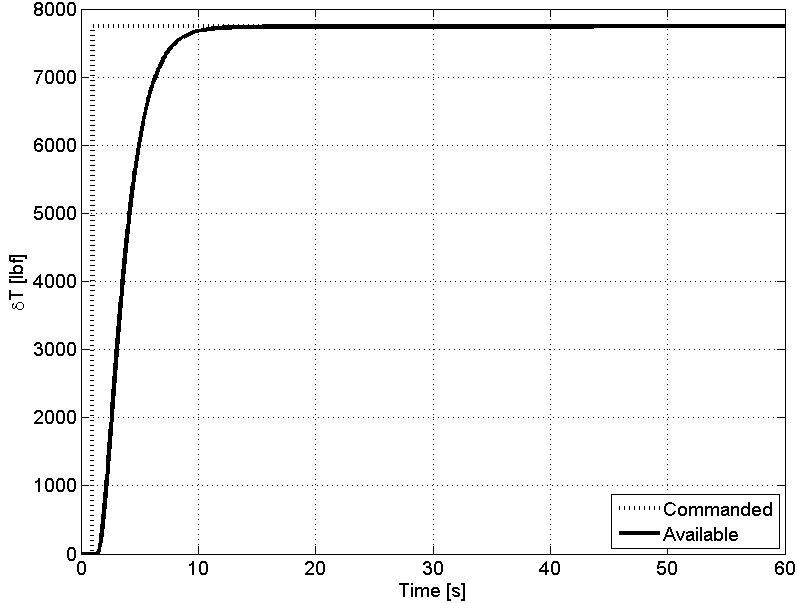}
 \caption{Commanded vs. available differential thrust}
 \label{fig:dT}
 \end{figure}

Unlike the rudder, due to delayed engine dynamics with time constant, there is a major difference in the commanded differential thrust and the available differential thrust as shown in Fig. \ref{fig:dT}.

It can be seen from Fig. \ref{fig:dT} that compared to the commanded differential thrust, the available differential thrust is equal in amount but longer in the time delivery. For a one degree step input on the rudder, the corresponding equivalent commanded and available differential thrust are 7737 lbf,  which is deliverable in ten second duration. Unlike the instantaneous control of the rudder input, there is a lag associated with the use of differential thrust as a control input. This is due to the lag in engine fluid transport and the inertia of the mechanical systems such as rotors and turbo-machinery blades \cite{NguyenStepanyan10}. This is a major constraint of plant dynamics and will be taken into account during the robust control system design phase in the following sections. 

\section{Open Loop System Response Analysis}\label{Open Loop}

\begin{figure}[htbp]
 \centering
 \includegraphics[width=5truein]{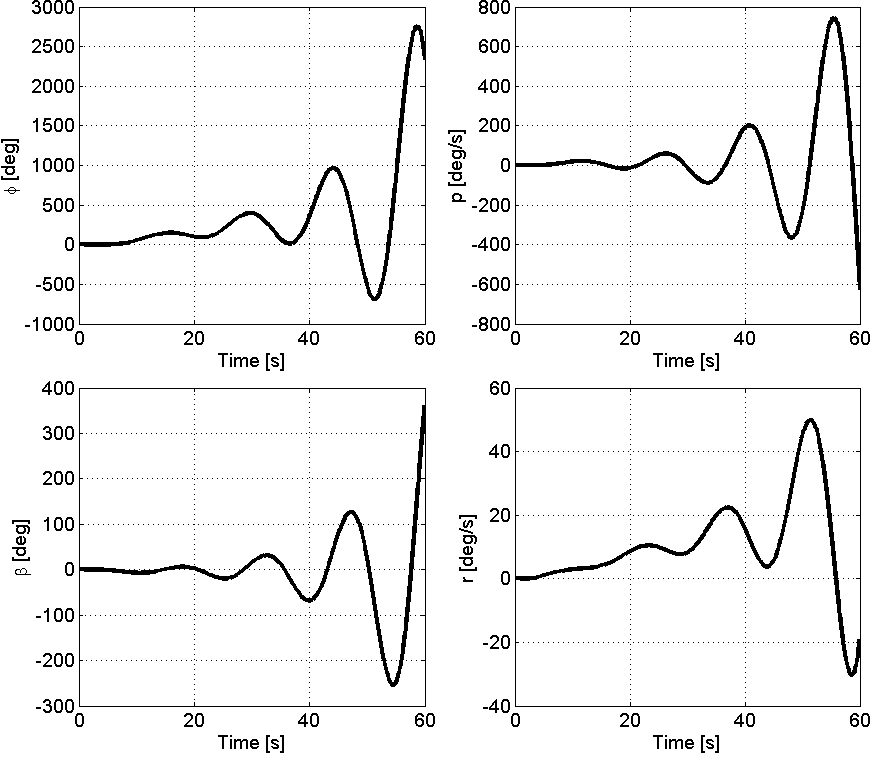}
 \caption{Open loop system response of the damaged aircraft}
 \label{fig:open_loop_response}
 \end{figure}

Following to this, the open loop response characteristics of the aircraft with a damaged vertical stabilizer to one degree step inputs from the ailerons and differential thrust are presented in Fig. \ref{fig:open_loop_response}. Here, it can be clearly seen that when the aircraft is damaged and the vertical stabilizer is lost, the aircraft response to the pilot's inputs is absolutely unstable in all four states (as it was also obvious (and expected) from the pole locations). This shows the pilot will not have much chance to stabilize the aircraft in time, which calls for a novel approach to save the damaged aircraft. This is another point where the second novel contribution of the paper is introduced: automatic control strategy to stabilize the aircraft,  which allows safe (i.e. intact) landing of the aircraft.

%
%
%
%
%
%
%
%
%

\section{Robust Control System Design based on $H_{\infty}$ Loop-Shaping Approach} \label{Robust}

\subsection{$H_{\infty}$ Loop-Shaping Design}

As having been witnessed in Section \ref{Open Loop}, the open loop responses of the damaged aircraft are unstable in all four lateral/directional states. This means the pilot will not have much chance to save the aircraft, which calls for a novel approach to save the damaged aircraft and to provide a safe landing. Therefore, in this section, the robust control system design based on $H_{\infty}$ loop-shaping approach is chosen as a means to stabilize the damaged aircraft due to its ability to suppress external disturbances and noise, and overall system robustness.

\begin{figure}[htbp]
 \centering
 \includegraphics[width=5truein]{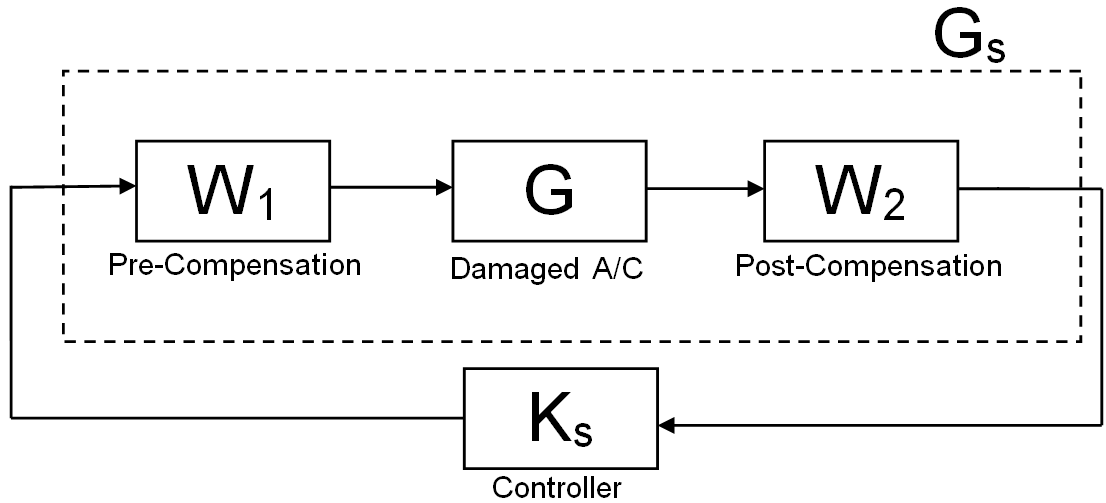}
 \caption{The damaged aircraft's shaped plant and controller}
 \label{fig:shaped_plant}
 \end{figure}

The damaged aircraft plant's open loop singular values are shaped by the pre- and post-compensation \cite{McFarlane_Glover90} as illustrated in Fig. \ref{fig:shaped_plant}.

As it can be seen from Fig. \ref{fig:shaped_plant}, $G$ is the open loop plant of the damaged aircraft while $W_1$ and $W_2$ are the pre- and post-compensation, respectively. The shaped plant is, therefore,

\begin{equation} \label{eq:G_s}
G_s=W_2GW_1
\end{equation}

In addition, the controller $(K_s)$ is synthesized by solving the robust stabilization problem for the shaped plant $(G_s)$ with a normalized left coprime factorization of $G_s={M_s}^{-1}N_s$, and the feedback controller for plant $G$ is, therefore, $K=W_1K_sW_2$ \cite{Skogestad_Postlethwaite96}.

\begin{figure}[htbp]
 \centering
 \includegraphics[width=6truein]{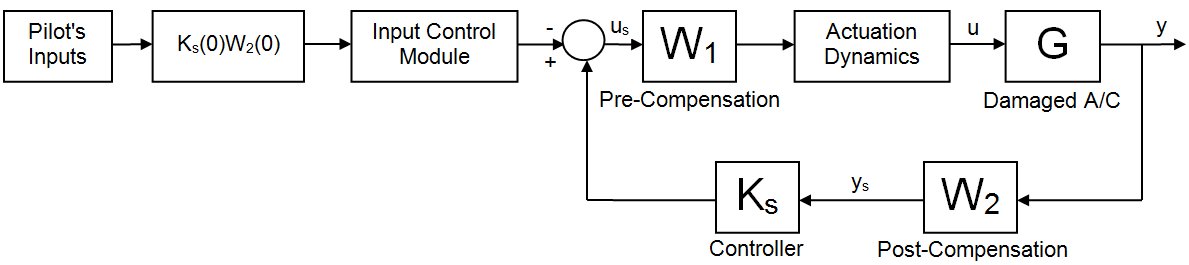}
 \caption{$H_{\infty}$ loop-shaping diagram}
 \label{fig:hinf_loop}
 \end{figure}

Next, the implemented $H_{\infty}$ loop-shaping diagram is shown in Fig. \ref{fig:hinf_loop}.

Figure \ref{fig:hinf_loop} provides an insight into the robust controller design based on the $H_{\infty}$ loop-shaping approach. The main goal of this controller is to stabilize all four lateral/directional states of the damaged aircraft, which are roll rate $(p)$, roll angle $(\phi)$, side-slip angle $(\beta)$, and yaw rate $(r)$. The control inputs of the damaged aircraft are ailerons $(\delta_a)$ and differential thrust $(\delta T)$. As it is obvious from Fig. \ref{fig:hinf_loop}, the pilot's inputs, which are one degree step inputs for both ailerons and rudder, will go through the pre-filter gain $K_s(0)W_2(0)$ and then the \emph{Input Control Module}, where the aileron control signal is routed through the saturation limits of $\pm$ 26 degrees \cite{faa_747} and the rudder control input is routed through the \emph{Differential Thrust Control Module}, where the rudder input is converted to differential thrust input as discussed previously in Section \ref{Differential Thrust}. It is also worth noting that there are also a differential thrust saturation limit of 43,729 lbf and a thrust generation rate limiter of 12,726 lbf/s imposed on the differential thrust control as discussed in Section \ref{Differential Thrust}. This is to make sure the control inputs are within the limits of both the ailerons and differential thrust. Furthermore, in order for the control efforts to be feasible in a real-life scenario, the control effort signals are also routed through the \emph{Actuation Dynamics}, where the same saturation and rate limits are imposed on the ailerons and differential thrust as discussed previously. 

The next step is to select the weighting functions, which should be taken into careful consideration due to the dimensions of the system. 
%

After an iterative process, the selected weighting functions for $W_1$ are as

\begin{equation}
{W_1}_{(11)}=\frac{4s + 1}{4s + 10}
\end{equation} 

\begin{equation}
{W_1}_{(22)}=\frac{50s + 5}{18s + 25}
\end{equation}  
  
\noindent We can then construct the \emph{system matrix} $W_1$ as

\begin{equation}
W_1= \left[
\begin{array}{ccccc}
-2.5000 & 0 & 2.5981 & 0 & 2.0000\\
0 & -1.3889 & 0 & 1.8922 & 0\\
-2.5981 & 0 & 3.0000 & 0 & 0\\
0 & -1.8922 & 0 & 2.7778 & 0\\
0 & 0 & 0 & 0 & -Inf
\end{array}\right]
\end{equation}  
  
The selected weighting functions for $W_2$ are as

\begin{equation}
{W_2}_{(11)}=\frac{16}{s + 16}
\end{equation} 

\begin{equation}
{W_2}_{(22)}={W_2}_{(33)}={W_2}_{(44)}=\frac{120}{s + 120}
\end{equation}

where the \emph{system matrix} of $W_2$ becomes

\begin{equation}
W_2= \left[
\begin{array}{ccccccccc}
-16.0000    &     0     &    0    &     0  &  4.0000   &      0   &      0    &     0  &  4.0000\\
0 & -120.0000    &     0    &     0   &      0  &  10.9545   &      0   &      0   &      0\\
0   &      0  & -120.0000   &      0   &      0   &      0 &   10.9545  &       0  &       0\\
0   &      0   &      0 & -120.0000   &      0   &      0   &      0  & 10.9545   &      0\\
4.0000   &      0   &      0    &     0   &      0   &      0   &      0    &     0    &     0\\
0 &  10.9545  &       0   &      0   &      0   &      0   &      0   &      0   &      0\\
0   &      0 &  10.9545   &      0   &      0   &      0   &      0   &      0   &      0\\
0   &      0  &       0 &  10.9545   &      0   &      0   &      0     &    0   &      0\\
0    &     0  &       0   &      0    &     0   &      0   &      0    &     0   &   -Inf
\end{array}\right]
\end{equation}  
  
As it is customary in general loop-shaping formulation, the maximum stability margin $(e_{max})$ is defined as a performance criterion
\begin{equation}
e_{max}=1/\gamma_{min}
\end{equation}  

\noindent where $\gamma_{min}$ is the $H_{\infty}$ optimal cost. For our analysis, the maximum stability margin is then obtained as

\begin{equation}
e_{max}= 0.2763
\end{equation}  

\noindent which is within the suggested optimal value bounds: $0.25<e_{max}<0.30$.

\begin{figure}[htbp] 
   \begin{minipage}[b]{0.5\linewidth}
    \centering
    \includegraphics[width=1\linewidth]{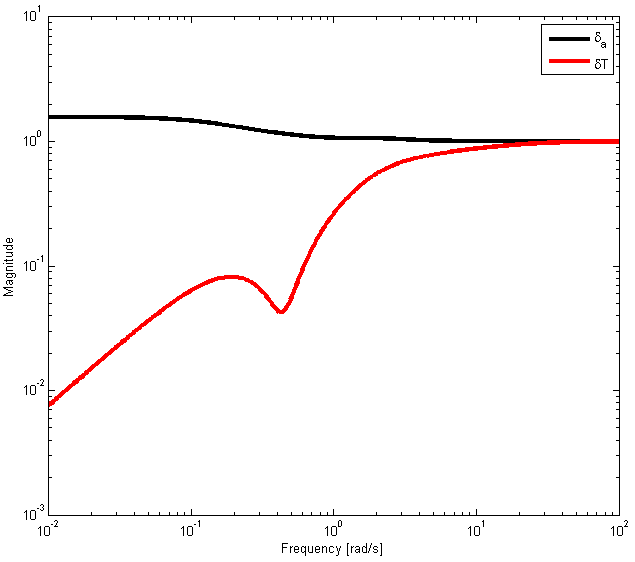}\\
    (a) Input sensitivity functions
    \vspace{2ex}
  \end{minipage}
  \begin{minipage}[b]{0.5\linewidth}
    \centering
    \includegraphics[width=1.06\linewidth]{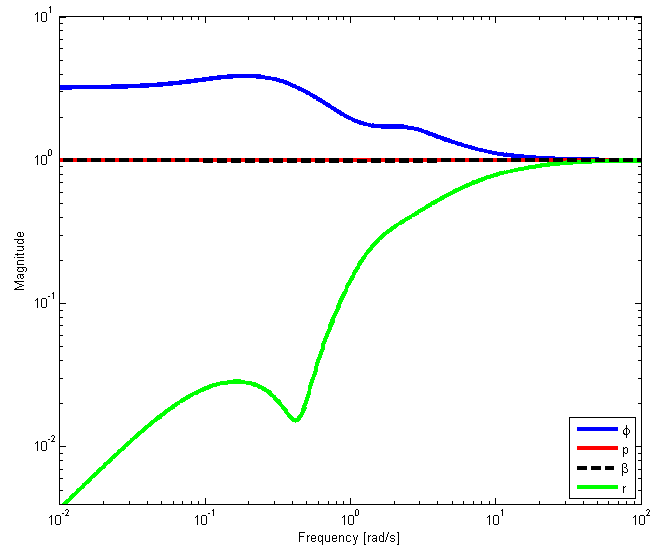}\\ 
    (b) Output sensitivity functions
    \vspace{2ex}
  \end{minipage} 
  \caption{Input and output sensitivity functions}
  \label{fig:Si_So}
 \end{figure}

Next, the Sensitivity (S) and Co-sensitivity (T) plots are constructed. 
Figure \ref{fig:Si_So} illustrates the input and output sensitivity functions.

%
 
\begin{figure}[htbp] 
   \begin{minipage}[b]{0.5\linewidth}
    \centering
    \includegraphics[width=1\linewidth]{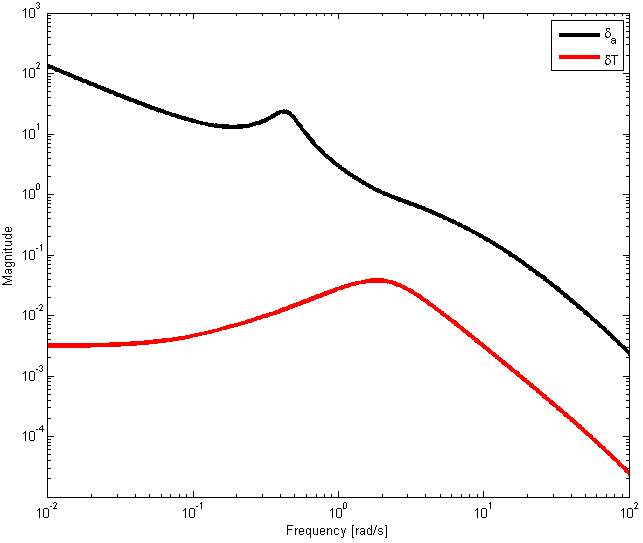}\\
    (a) Input responses of shaped plant
    \vspace{2ex}
  \end{minipage}
  \begin{minipage}[b]{0.5\linewidth}
    \centering
    \includegraphics[width=1.03\linewidth]{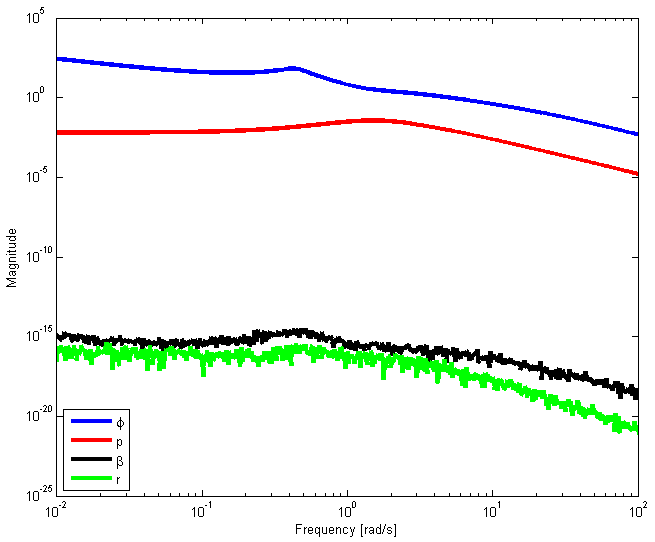}\\ 
    (b) Output responses of shaped plant
    \vspace{2ex}
  \end{minipage} 
  \caption{Input and output responses of shaped plant}
  \label{fig:Li_Lo}
 \end{figure}

Next, we can investigate the shaped plant behaviors for inputs and outputs in Fig. \ref{fig:Li_Lo}.

%
 
\noindent From Fig. \ref{fig:Li_Lo}, we can see that for the aileron control input $(\delta_a)$, the shaped plant has very good response characteristics of having high gain at low frequencies for tracking and low gain at high frequencies for disturbance rejection. For the differential thrust control input $(\delta T)$, although the gain at low frequencies is not as high as that of the aileron control input, but the gain is relatively linear, which makes the differential thrust control quite predictable for the pilots. At high frequencies, the gain roll-off is also linear and quite similar to that of the aileron control input, which is helpful at rejecting disturbances. In addition, Fig. \ref{fig:Li_Lo} also shows the output responses of the shaped plant in all four lateral/directional states, from which we can see that the implemented controller based on $H_{\infty}$ loop-shaping approach can augment the system in two groups: roll angle and roll rate ($\phi$ and $p$) and side-slip angle and yaw rate ($\beta$ and $r$). It is also obvious that roll angle and roll rate ($\phi$ and $p$) have higher gains than side-slip angle and yaw rate ($\beta$ and $r$), which is expected when the aircraft loses its vertical stabilizer. 

\begin{figure}[htbp]
 \centering
 \includegraphics[width=4truein]{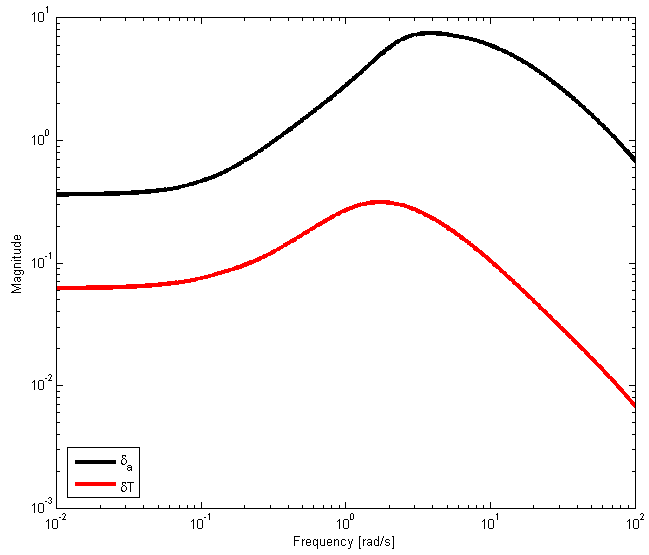}
 \caption{Augmented control action by the $H_{\infty}$ loop-shaping controller}
 \label{fig:K_aug}
 \end{figure}

\begin{figure}[htbp] 
   \begin{minipage}[b]{0.5\linewidth}
    \centering
    \includegraphics[width=1.03\linewidth]{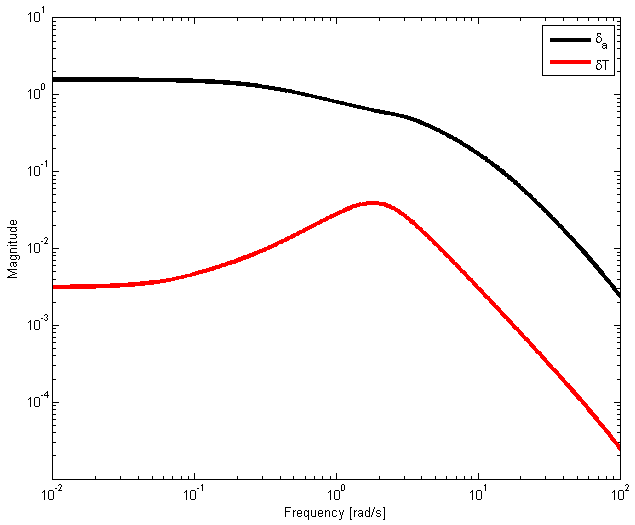}\\
    (a) Co-sensitivity plot for input responses of shaped plant
    \vspace{2ex}
  \end{minipage}
  \begin{minipage}[b]{0.5\linewidth}
    \centering
    \includegraphics[width=1.02\linewidth]{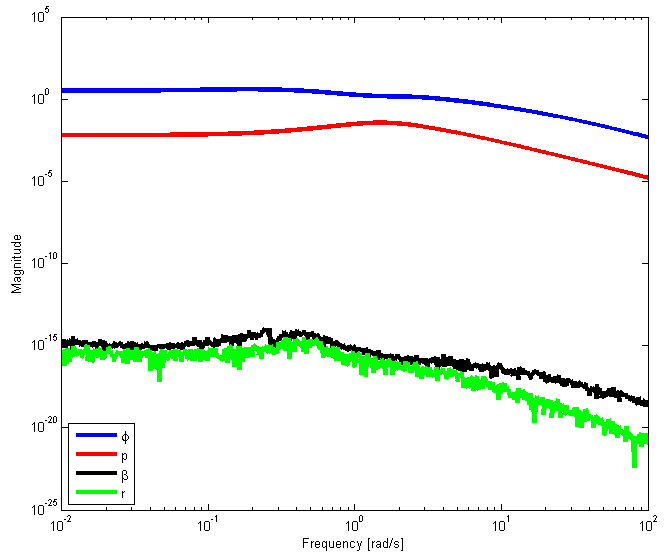}\\ 
    (b) Co-sensitivity plot for output responses of shaped plant
    \vspace{2ex}
  \end{minipage} 
  \caption{Co-sensitivity plot for input and output responses of shaped plant}
  \label{fig:Ti_To}
 \end{figure}

In addition, it is possible to look at the augmented control action by the $H_{\infty}$ loop-shaping controller in Fig. \ref{fig:K_aug}. Furthermore, the Co-sensitivity plots related to the shaped plant behaviors for the inputs and outputs are presented in Fig. \ref{fig:Ti_To}. Fig \ref{fig:K_aug} and Fig. \ref{fig:Ti_To} prove that the controller based on $H_{\infty}$ loop-shaping approach is able to achieve desirable closed loop response characteristics in the frequency domain in terms of robustness and stability.

%
%
 
\begin{figure}[htbp] 
   \begin{minipage}[b]{0.5\linewidth}
    \centering
    \includegraphics[width=1\linewidth]{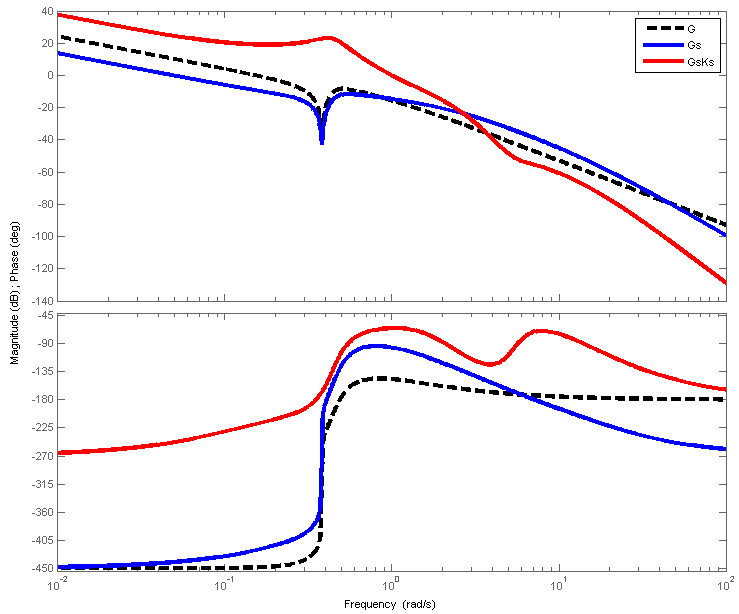}\\ 
    (a) Transfer function $\frac{\phi}{\delta_a}$
    \vspace{2ex}
  \end{minipage}
  \begin{minipage}[b]{0.5\linewidth}
    \centering
    \includegraphics[width=1.01\linewidth]{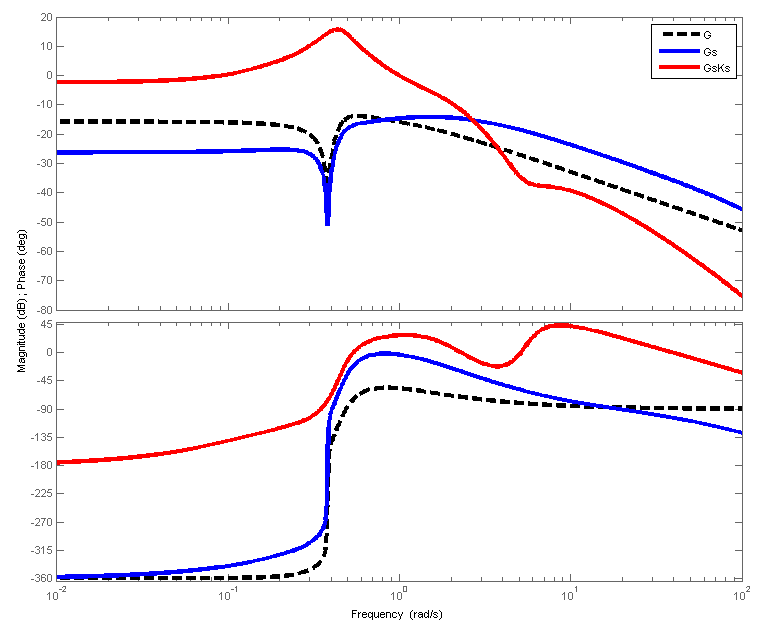} \\
    (b) Transfer function $\frac{p}{\delta_a}$
    \vspace{2ex}
  \end{minipage} 
  \begin{minipage}[b]{0.5\linewidth}
    \centering
    \includegraphics[width=1.008\linewidth]{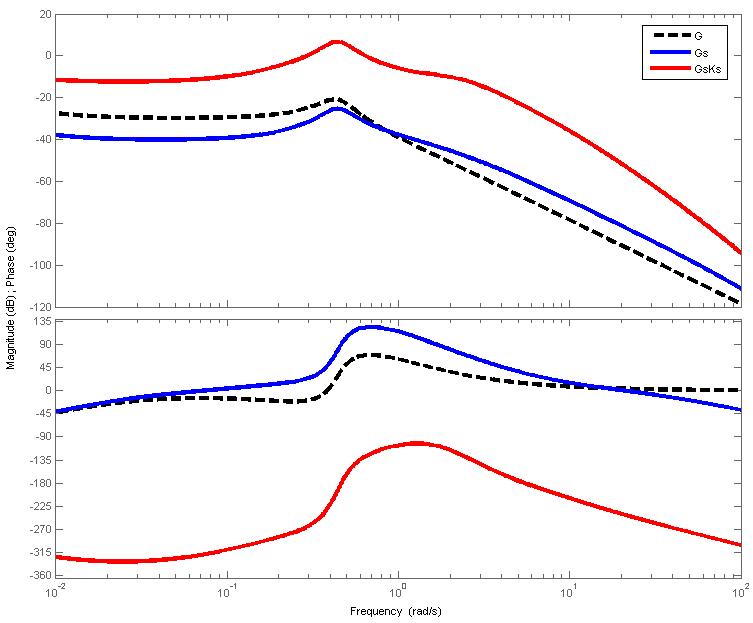}\\ 
    (c) Transfer function $\frac{\beta}{\delta_a}$
    \vspace{2ex}
  \end{minipage}
  \begin{minipage}[b]{0.5\linewidth}
    \centering
    \includegraphics[width=1\linewidth]{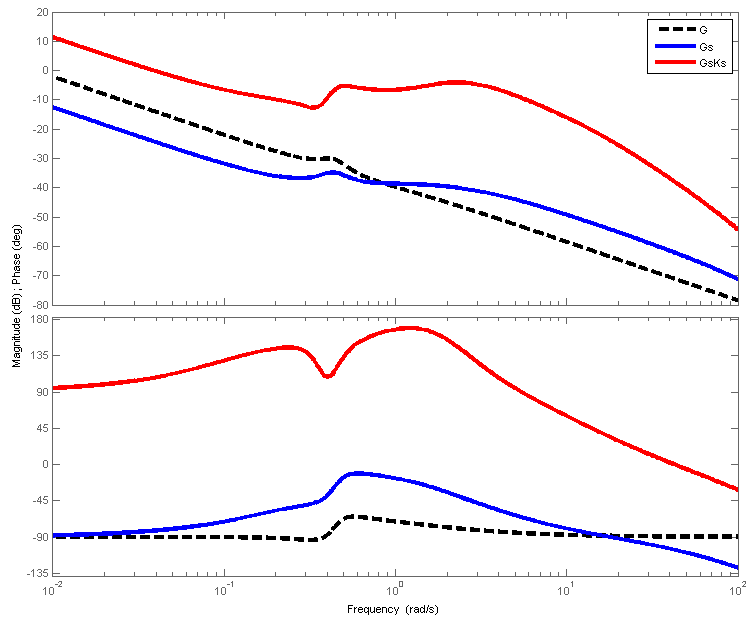}\\ 
    (d) Transfer function $\frac{r}{\delta_a}$
    \vspace{2ex}
  \end{minipage} 
  \caption{Open-loop $(G)$ vs. shaped $(G_s)$ vs. robustified $(G_sK_s)$ plant for aileron input}
  \label{fig:loop_da}
\end{figure}

\begin{figure}[htbp] 
    \begin{minipage}[b]{0.5\linewidth}
    \centering
    \includegraphics[width=1\linewidth]{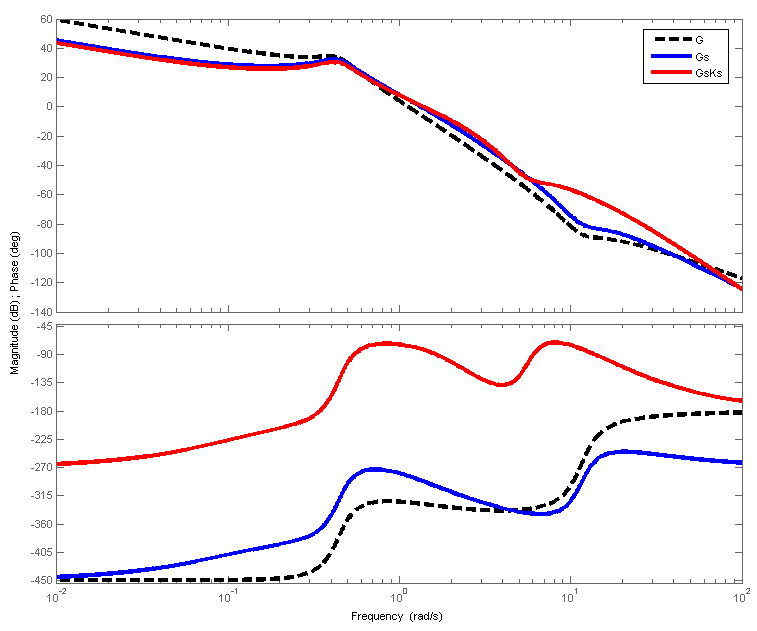}\\
    (a) Transfer function $\frac{\phi}{\delta T}$
    \vspace{2ex}
  \end{minipage}
  \begin{minipage}[b]{0.5\linewidth}
    \centering
    \includegraphics[width=1\linewidth]{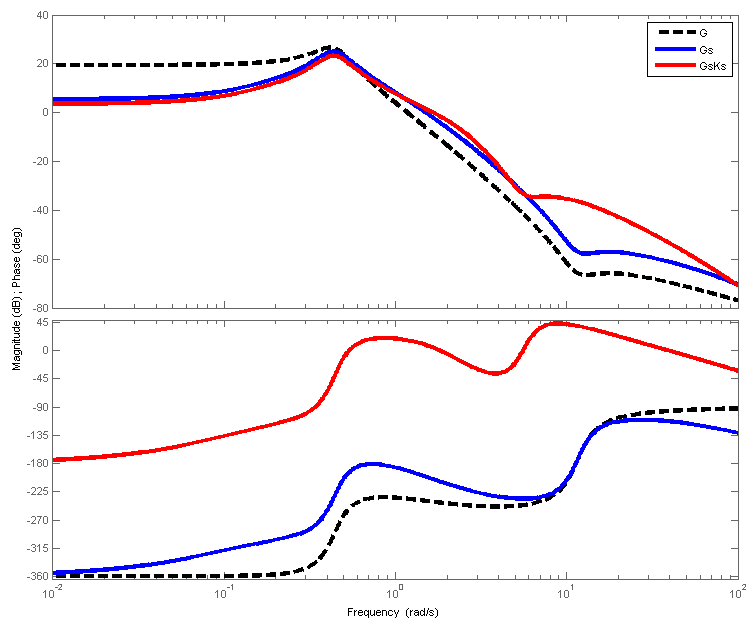}\\ 
    (b) Transfer function $\frac{p}{\delta T}$
    \vspace{2ex}
  \end{minipage} 
  \begin{minipage}[b]{0.5\linewidth}
    \centering
    \includegraphics[width=1.01\linewidth]{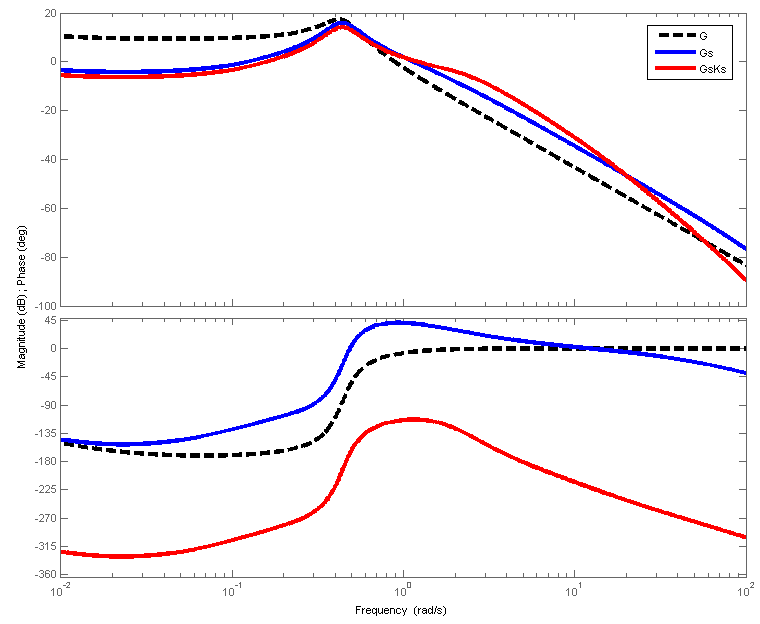}\\ 
    (c) Transfer function $\frac{\beta}{\delta T}$
    \vspace{2ex}
  \end{minipage}
  \begin{minipage}[b]{0.5\linewidth}
    \centering
    \includegraphics[width=1\linewidth]{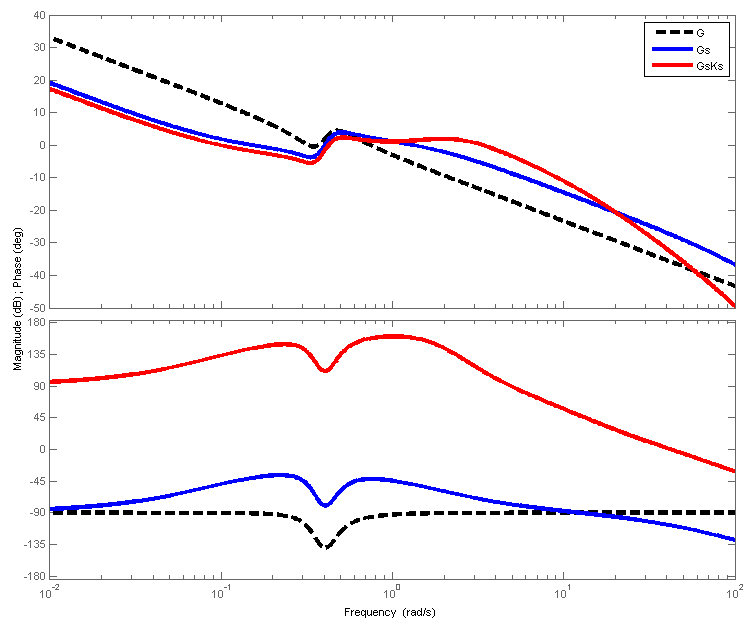} \\
    (d) Transfer function $\frac{r}{\delta T}$ 
    \vspace{2ex}
  \end{minipage} 
    \caption{Open-loop $(G)$ vs. shaped $(G_s)$ vs. robustified $(G_sK_s)$ plant for differential thrust input}
  \label{fig:loop_dT}
\end{figure}

Next, the comparison of the open-loop, shaped, and robustified plant response is carried out to investigate the effect of the $H_{\infty}$ loop-shaping on the plant dynamics as illustrated in Fig.  \ref{fig:loop_da} and Fig. \ref{fig:loop_dT}. 

As seen in Fig.  \ref{fig:loop_da} and Fig. \ref{fig:loop_dT}, compared to the open-loop plant response, the performance of the damaged aircraft is further improved by controller. 

\subsection{Allowable Gain/Phase Variations: Disk Margin Analysis}

\begin{figure}[htbp]
 \centering
 \includegraphics[width=4truein]{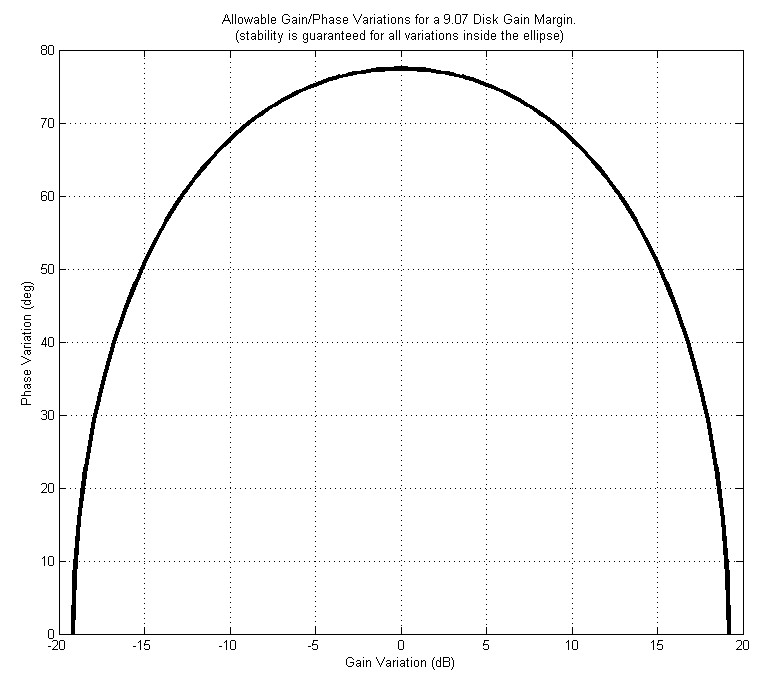}
 \caption{Disk gain and phase margins for inputs}
 \label{fig:diskmargin_input}
 \end{figure}

\begin{figure}[htbp]
 \centering
 \includegraphics[width=4truein]{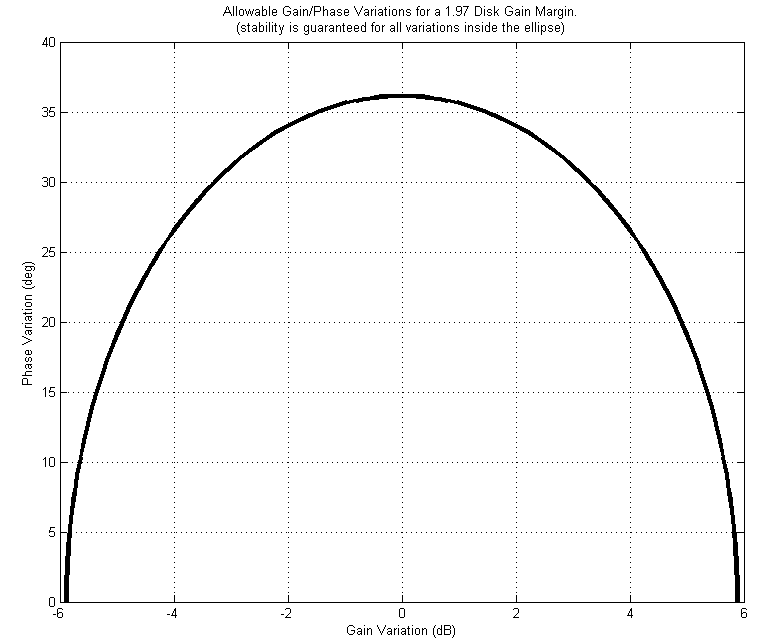}
 \caption{Disk gain and phase margins for outputs}
 \label{fig:diskmargin_output}
 \end{figure}

Furthermore, it is also worth checking the allowable gain and phase margin variations (i.e. associated uncertainty balls in gain and phase margins for a robust response and guaranteed stability). Obtained results are described in Fig. \ref{fig:diskmargin_input} and Fig. \ref{fig:diskmargin_output}. 

From Fig. \ref{fig:diskmargin_input} and Fig. \ref{fig:diskmargin_output}, it can be clearly seen that with the $H_{\infty}$ loop-shaping control system design, both the inputs and outputs can achieve desirable and safe stability margins. 

\subsection{Simulation Results}
 
\begin{figure}[htbp]
 \centering
 \includegraphics[width=6truein]{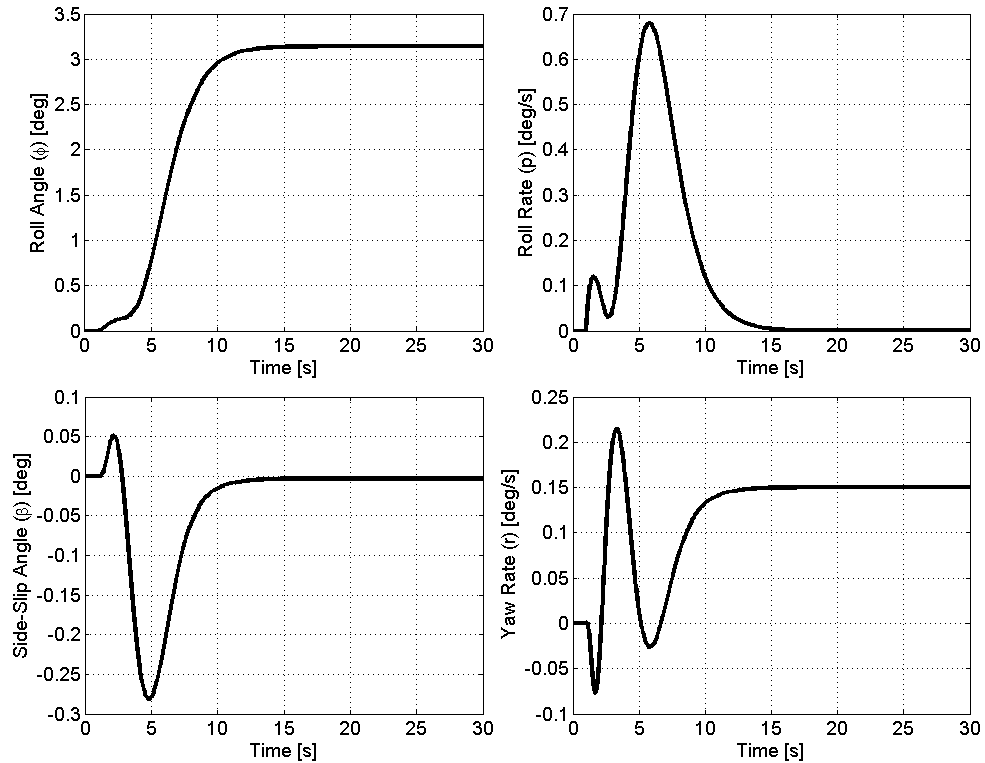}
 \caption{Closed loop time domain responses of the damaged aircraft}
 \label{fig:hinf_state}
 \end{figure}

\begin{figure}[htbp]
 \centering
 \includegraphics[width=4truein]{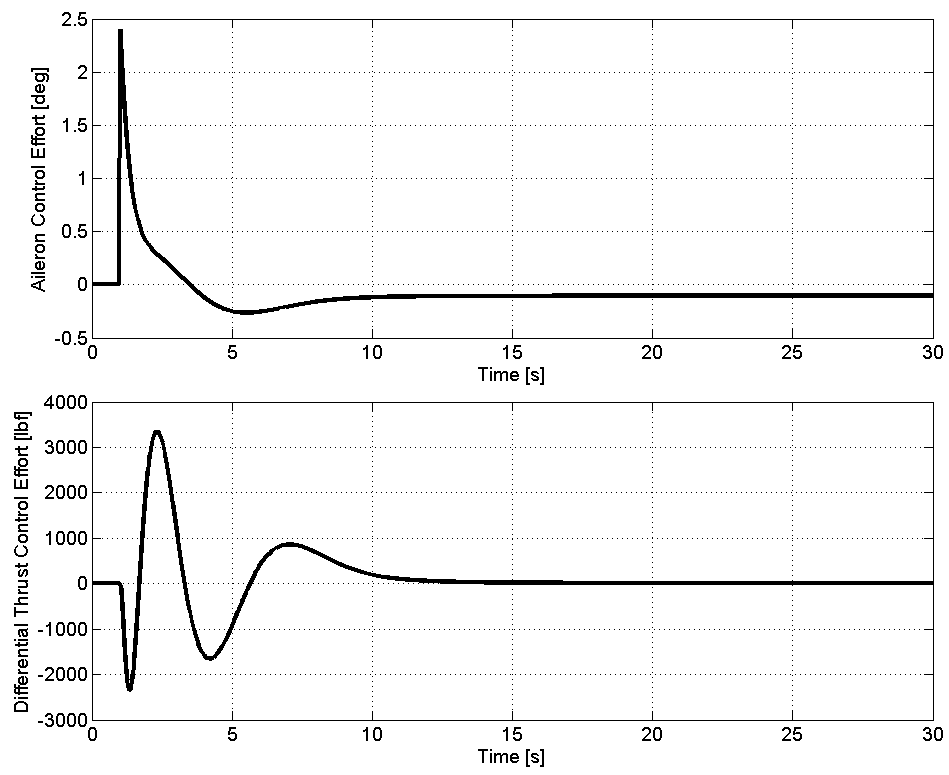}
 \caption{Control efforts}
 \label{fig:hinf_control_efforts}
 \end{figure}

Now it is time to look at the time domain behaviors of our system, which are shown in Fig. \ref{fig:hinf_state}. It is worth mentioning that the control inputs for both plants are one degree step inputs for both the ailerons and differential thrust to simulate an extreme scenario test to see whether the damaged aircraft utilizing differential thrust can hold itself in a continuous yawing and banking maneuver without becoming unstable and losing control. From Fig. \ref{fig:hinf_state}, it is obvious that after only 15 seconds, all four states of the aircraft's lateral/directional dynamics reach steady state values, which means that the controller can stabilize the damaged aircraft in only 15 seconds. However, this comes at the cost of the control efforts as shown in Fig. \ref{fig:hinf_control_efforts}, which are still under the control limits and without any saturation of the actuators.

As discussed previously in this section, in order to have a feasible control strategy in real-life situation, limiting factors are imposed on the aileron and differential thrust control efforts. The aileron deflection is limited at $\pm$26 degrees \cite{faa_747}. For differential thrust, a differential thrust saturation is set at 43,729 lbf, which is the difference of the maximum thrust and trimmed thrust values of the JT9D-7A engine. In addition, a rate limiter is also imposed on the thrust response characteristic at 12,726 lbf/s as discussed in Section \ref{Differential Thrust}.

The aileron control effort, as indicated by Fig. \ref{fig:hinf_control_efforts}, calls for the maximum deflection of approximately 2.4 degrees and reaches steady state at approximately -0.1 degree of deflection after 15 seconds. This aileron control effort is very reasonable and achievable if the ailerons are assumed to have instantaneous response characteristics by neglecting the lag from actuators or hydraulic systems. The differential thrust control effort demands a maximum differential thrust of approximately 3350 lbf, which is well within the thrust capability of the JT9D-7A engine, and the differential thrust control effort reaches steady state at around 15 lbf after 15 seconds. Therefore, it can be concluded that the robust control system design based on the $H_{\infty}$ loop-shaping approach is proven to be able to stabilize and save the damaged aircraft.

\section{Robustness and Uncertainty Analysis} \label{Robustness}

\begin{figure}[htbp!]
 \centering
 \includegraphics[scale=0.5]{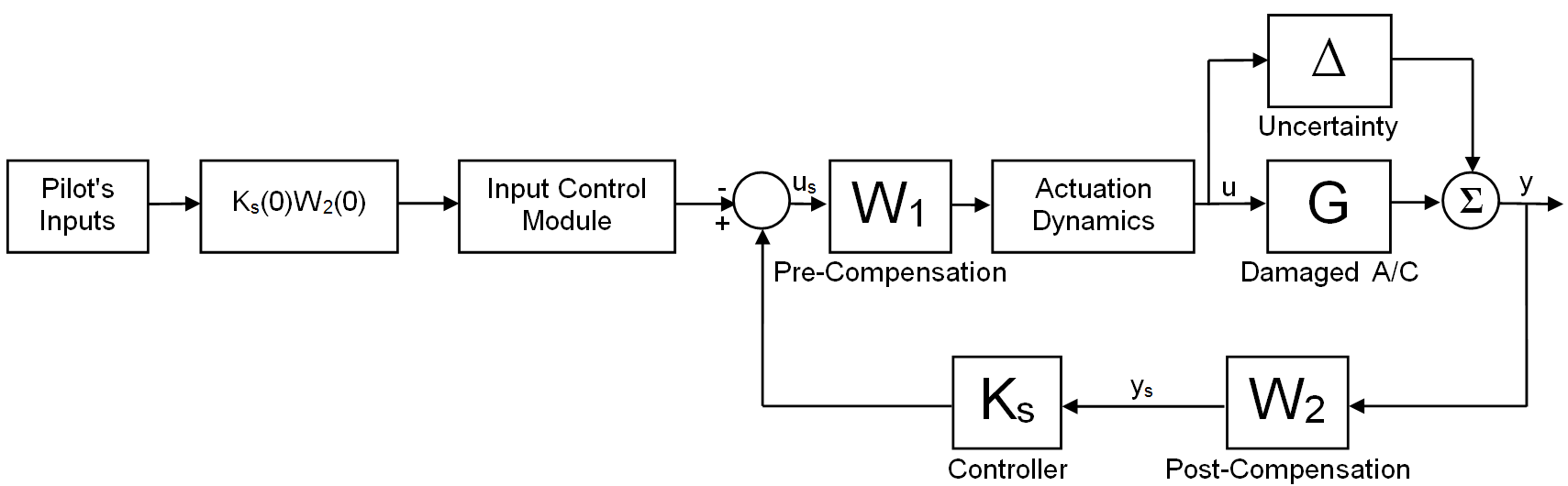}
 \caption{$H_{\infty}$ loop-shaping diagram in the presence of uncertainty}
 \label{fig:robust_hinf_loop}
 \end{figure}

The robustness of the $H_{\infty}$ loop-shaping robust differential thrust control system design presented in this paper is investigated by the introduction of 30\% of full block, additive uncertainty into the plant dynamics of the damaged aircraft to test the performance of the aircraft in the presence of uncertainty. Fig. \ref{fig:robust_hinf_loop} shows the logic behind the robust control system design in the presence of uncertainty.

\begin{figure}[htbp!]
 \centering
 \includegraphics[width=6truein]{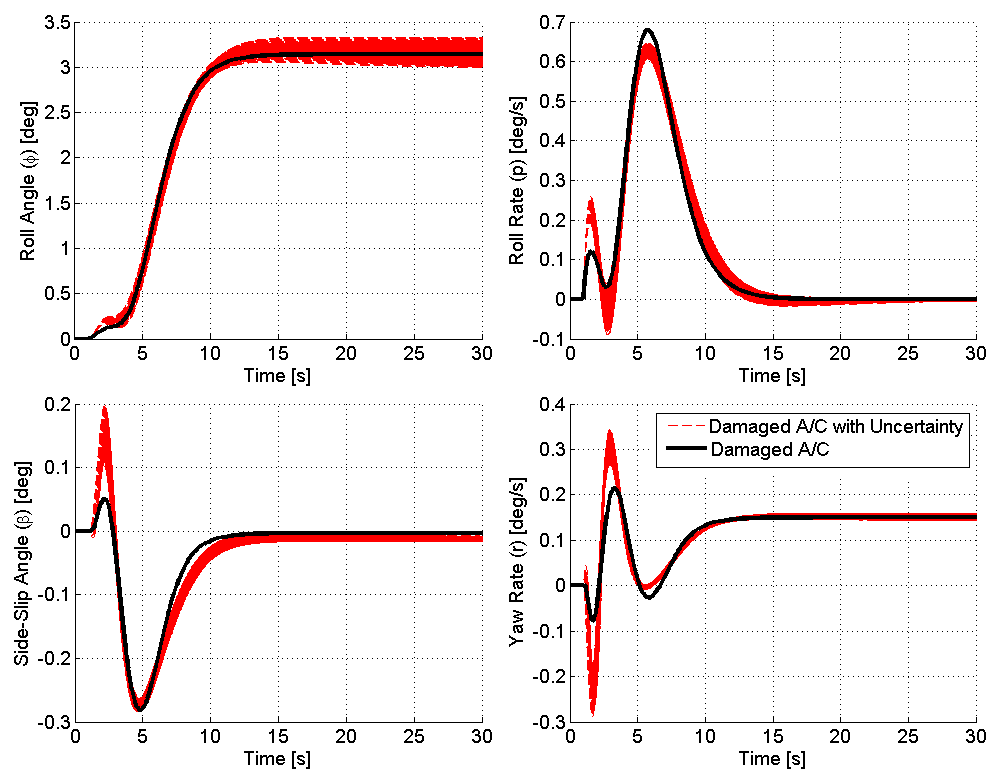}
 \caption{Closed loop time domain responses of the damaged aircraft in the presence of 30\% uncertainty}
 \label{fig:robust_output_U_30}
 \end{figure}

One thousand Monte-Carlo simulations were conducted to test the robustness of the damaged plant in the presence of uncertainty. The state responses in the presence of 30\% of uncertainty are shown in Fig. \ref{fig:robust_output_U_30}, where it is obvious that the robust control system design is able to perform well under given uncertain conditions and the damaged aircraft has stable steady state responses within only 15 seconds. In that sense, the uncertain plant dynamics are well within the expected bounds.
 
\begin{figure}[htbp!]
 \centering
 \includegraphics[width=4truein]{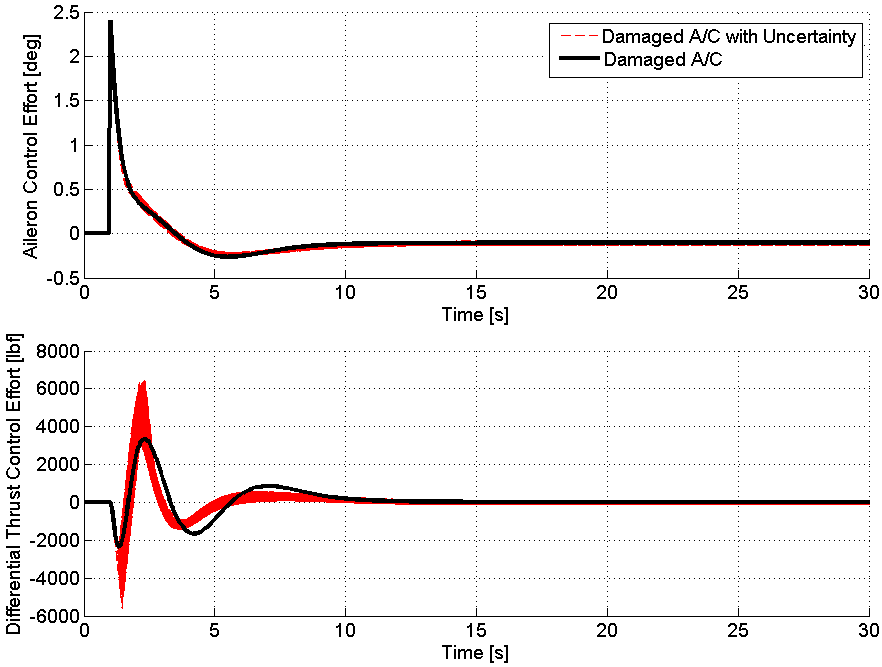}
 \caption{Control efforts in the presence of 30\% uncertainty}
 \label{fig:robust_control_effort_U_30}
 \end{figure}

However, these favorable characteristics come at the expense of the control effort from the ailerons and differential thrust as shown in Fig. \ref{fig:robust_control_effort_U_30}.

According to Fig. \ref{fig:robust_control_effort_U_30}, when there is 30\% full block, additive uncertainty, the aileron control demands the maximum deflection of approximately 2.4 degrees and reaches steady state between -0.12 and -0.09 degree after 15 seconds, which is quite similar to the required aileron control effort when there is no uncertainty. Therefore, the aileron control effort demands are reasonable and feasible due to the limiting factor of $\pm$ 26 degrees of the aileron deflection \cite{faa_747} and the assumption that ailerons have instantaneous response characteristics by neglecting the lag from actuators or hydraulic systems. 

As for differential thrust, when there is 30\% uncertainty, the differential thrust control demands at maximum approximately 6500 lbf, which is within the thrust capability of the JT9D-7A engine, and the differential thrust control effort reaches steady state between -100 lbf (negative differential thrust means $T1<T4$) and 100 lbf (positive differential thrust means $T1>T4$) after 15 seconds. Compared to the case when there is no uncertainty, the demanded differential thrust associated with uncertain plant dynamics is higher in both magnitude and rate. It is also obvious from Fig. \ref{fig:robust_control_effort_U_30} that in a few cases, the differential thrust control effort demand hit the thrust generation rate limiter, which is set at 12,726 lbf/s for the JT9D-7A engine, but fortunately, the control system is so robust that throughout 1000 Monte-Carlo simulations, it can stabilize the aircraft's uncertain plant dynamics. Again, due to the differential thrust saturation set at 43,729 lbf and the thrust response limiter set at 12,726 lbf/s, this control effort of differential thrust in the presence of uncertainty is achievable in real life situation.

\section{Conclusion and Future Work} \label{Conclusion}
This paper studied the utilization of differential thrust as a control input to help a Boeing 747-100 aircraft with a damaged vertical stabilizer regain its lateral/directional stability. The motivation of this research study is to improve the safety of air travel in the event of losing the vertical stabilizer by providing control means to safely control and/or land the aircraft.

Throughout this paper, the necessary damaged aircraft model was constructed, where lateral/directional equations of motion were revisited to incorporate differential thrust as a control input for the damaged aircraft. Then the open loop plant dynamics of the damaged aircraft was investigated. The engine dynamics of the aircraft was modeled as a system of differential equations with engine time constant and time delay terms to study the engine response time with respect to a commanded thrust input. Next, the novel differential thrust control module was presented to map the rudder input to differential thrust input. The $H_{\infty}$ loop-shaping approach based robust control system design's ability to stabilize the damaged aircraft was proven as investigation results demonstrated that the damaged aircraft was able to reach steady state stability within only 15 seconds under feasible control efforts. Further analysis on robustness showed that the $H_{\infty}$ loop-shaping approach based robust differential thrust methodology was also able to stabilize the uncertain plant dynamics in the presence of 30\% full block, additive uncertainty associated with the damaged aircraft dynamics. 

Through listed analyses above, the ability to save the damaged aircraft by the robust differential thrust control strategy  has been demonstrated in this paper. This framework provides an automatic control methodology to save the damaged aircraft and avoid the dangerous coupling of the aircraft and pilots, which led to crashes in a great number of commercial airline incidents. Furthermore, it has also been concluded that due to the heavy dependence of the differential thrust generation on the engine response, in order to better incorporate the differential thrust as an effective control input in a life-saving scenario, major developments in engine response characteristics are also desired to better assist such algorithms.

Future work includes the implementation of such control algorithms into a small scaled unmanned aerial vehicle (UAV) for real-life testing.


\end{document}